\documentclass[12pt]{article}

\usepackage{indentfirst}
\usepackage{footmisc}

\usepackage[latin1]{inputenc}
\usepackage{amsmath,amssymb}
\usepackage{amsmath,accents}

\usepackage{alltt,graphics,graphpap,color}
\usepackage[dvips]{graphicx}

\usepackage{amsfonts}
\usepackage{amscd}

\usepackage{mathrsfs}
\usepackage{hyperref}
\usepackage{ marvosym }

\hypersetup{linktoc=page,colorlinks=true,linkcolor=red,  citecolor=blue, }

\newtheorem{Theorem}{Theorem}[section]

\newtheorem{Proposition}[Theorem]{Proposition}
\newtheorem{Lemma}[Theorem]{Lemma}
\newtheorem{Corollary}[Theorem]{Corollary}

\newtheorem{Example}[Theorem]{Example}
\newtheorem{Notation}[Theorem]{Notation}

\newcommand{\proof}{\noindent\emph{Proof. }}
\newcommand{\cvd}{\hfill$\square$ \bigskip}

\newcommand{\mm}{\mathcal {M}}

\newcommand{\kkk}{\mathbb K}
\newcommand{\rrr}{\mathbb R}
\newcommand{\ccc}{\mathbb C}
\newcommand{\hhh}{\mathbb H}
\newcommand{\sss}{\mathbb S}
\newcommand{\kh}{\mathbb {KH} }

\newcommand{\cth}{\frac{\cosh(\rho)}{\sinh(\rho)}}
\newcommand{\tgh}{\frac{\sinh(\rho)}{\cosh(\rho)}}
\newcommand{\sih}{\sinh(\rho)}
\newcommand{\coh}{\cosh(\rho)}
\newcommand{\coq}{\cosh^2(\rho)}
\newcommand{\siq}{\sinh^2(\rho)}

\newcommand{\pma}{\psi(m+a)}
\newcommand{\pmm}{\psi(m)}
\newcommand{\eo}{\frac{2}{\psi(m+a)}}
\newcommand{\deo}{\frac{4}{\psi(m+a)}}
\newcommand{\LLL}{\frac{\psi'}{\psi^2}\Delta}

\newcommand{\aao}{\accentset{\ \circ}{\left|A\right|}^2}

\newcommand{\gtf}{\accentset{\phantom{aaa.}\circ}{|\nabla A|}^2}
\newcommand{\gtfq}{\accentset{\phantom{aaa.}\circ}{|\nabla A|}^4}
\newcommand{\tf}{\accentset{\circ}{h}}

\begin{document}

\title{Nonhomogeneous expanding flows in hyperbolic spaces}

\author{\sc Giuseppe Pipoli }
\date{}

\maketitle

\vspace{1,5cm}

{\small \noindent {\bf Abstract:} A recent paper \cite{CGT} studies the evolution of star-shaped mean convex hypersurfaces of the Euclidean space by a class of nonhomogeneous expanding curvature flows. In the present paper we consider the same problem in the real, complex and quaternionic hyperbolic spaces, investigating how the richer geometry of the ambient space affects the evolution. In every case the initial conditions are preserved and the long time existence of the flow is proven. The geometry of the ambient space influences the asymptotic behaviour of the flow: after a suitable rescaling the induced metric converges to a conformal multiple of the standard Riemannian round metric of the sphere if the ambient manifold is the real hyperbolic space, otherwise it converges to a conformal multiple of the standard sub-Riemannian metric on the odd-dimensional sphere. Finally, in every cases, we are able to construct infinitely many examples such that the limit does not have constant scalar curvature.} 

\medskip

 \bigskip

\section{Introduction}

Recently Li Chen, Xi Guo and Qiang Tu \cite{CGT} considered the evolution of closed, mean convex and star-shaped hypersurfaces of the Euclidean space by a class of nonhomogeneous expanding curvature flows modelled on the inverse mean curvature flow. They proved that these flows have a common behaviour: the initial conditions are preserved, the flow is defined for any positive time and the induced metric converges, after rescaling, to the usual round metric on the sphere regardless of the initial datum.

In recent years many results about nonhomogeneous curvature flows in the Euclidean space appear. Different types of speed have been studied and many different problems addressed. Just to mention some of the most recent results, Sinestrari and his collaborators produced convexity estimates \cite{AS}, considered volume and area preserving flows \cite{BS} and ancient solutions for a very general class of expanding flows \cite{RS}. Moreover, McCoy considered contracting nonhomogeneous flows \cite{MC1} and their self-similar solutions \cite{MC2}.

Despite a growing interest in nonhomogeneous flows, the literature about their evolution in Riemannian manifolds is still at the beginning. For the best of our knowledge the first paper about nonhomogeneous flows in a curved space is that one of Bertini with the author of the present paper \cite{BP}, where we consider volume and area preserving flow in the real hyperbolic space. The ambient manifold is always a space form in all the recent papers on the subject: \cite{AW,ACW,LL,LZ}. 

The goal of this paper is to extend the results of \cite{CGT} studying analogous flows in the real, complex and quaternionic hyperbolic spaces and exploring how a richer geometry can affect the evolution. 

Let $\kkk$ be either the real field $\rrr$, or the complex field $\ccc$, or the algebra of quaternions $\hhh$. Let $\kh^n$ be the $\kkk$-hyperbolic space endowed with its standard Riemannian metric. Let $F:\mm\times [0,T)\rightarrow\kh^n$ be a one-parameter family of smooth embeddings such that $F_0(\cdot)=F(\cdot,0)$ is a given hypersurface and $F$ evolves by
\begin{equation}\label{def flow}
\frac{\partial F}{\partial t} = \frac{1}{\psi(H)}\nu,
\end{equation}
where $\nu$ is the unit outward normal vector of $F$, $H$ is its mean curvature and $\psi:[0,\infty)\rightarrow\rrr$ is a continuous function $C^2$ differentiable in $(0,\infty)$ which satisfies the following structural conditions:

\begin{eqnarray} 
\nonumber i) &&\psi(0)=0,\ \psi(x)>0,\ \psi'(x)>0,\qquad \forall x>0;\\
\label{condition} ii)&& \frac{x\psi'(x)}{\psi(x)}\leq 1,\qquad \forall x>0;\\
\nonumber iii)&& \psi''(x)\psi(x)-2(\psi'(x))^2\leq 0,\qquad \forall x>0.
\end{eqnarray}

These conditions are part of the properties of the speed considered in \cite{CGT}, but our class is more general. In fact, in our case, we will prove that $H$ cannot converges to zero, hence we do not need to prescribe the behaviour of $\psi'$ when $x$ is tends to $0$. As said in \cite{CGT} \eqref{condition} includes homogeneous functions, as the suitable powers of $x$, but also many nonhomogeneous ones, as $\ln(1+x)$, or $\sum_{i=1}^kc_ix^{-p_i}$ with $c_i>0$ and $0<p_i\leq 1$. Note that it could happen that $\psi$ is not concave, but Condition \eqref{condition} iii) is equivalent to the fact that the function $\frac{1}{\psi}$ is convex.

Before giving the precise statement of our main result we need to introduce some useful notations in order to consider the three cases at once. For any $\kkk$ we define:
\begin{eqnarray}\label{def a}
a \ :=\  \dim_{\rrr}\kkk-1& =& \left\{\begin{array}{ccc}
0 & \text{if}& \kkk=\rrr,\\
1& \text{if}& \kkk=\ccc,\\
3& \text{if}& \kkk=\hhh;
\end{array}\right.\\
\label{def m}
m \ := \ \dim_{\rrr}\kh^n-1 & =& (a+1)n-1.
\end{eqnarray}

Clearly $m$ is the dimension of any real hypersurface in $\kh^n$. Let $\sigma$ be the usual round metric on the sphere. If $\kkk\neq\rrr$ we need to introduce some other metrics on $\sss^m$. The action of $\sss^a$ on $\sss^m$ induces the Hopf fibration $\pi:\sss^m\rightarrow\mathbb{KP}^{n-1}$. Let $\mathcal H$ be the horizontal distribution of the Riemannian submersion $\pi$. We denote with $\sigma_{\kkk}$ the sub-Riemannian metric on $\sss^m$ which coincides with $\sigma$ on $\mathcal H\times\mathcal H$ (and it is not defined outside). We will call it the \emph{standard sub-Riemannian metric} on $\sss^m$. The main result of this paper is the following.

\begin{Theorem}\label{main}
Let $\mm_0$ be a closed star-shaped, mean convex hypersurface of $\kh^n$. If $\kkk=\rrr$ suppose that $n\geq 3$, otherwise suppose that $n\geq 2$ and $\mm_0$ is $\sss^a$-invariant. Let $\mm_t$ be the evolution of $\mm_0$ along the nonhomogeneous flow \eqref{def flow} , where $\psi$ satisfies Conditions \eqref{condition}. Let $g_t$ be the induced metric on $\mm_t$ and consider the rescaled metric
$$
\tilde g_t=|\mm_t|^{-\frac{2}{m+a}}g_t.
$$ 
Then:
\begin{itemize}
\item[(1)] $\mm_t$ is star-shaped, mean convex and $\sss^a$-invariant for any time the flow is defined;
\item[(2)] the flow is defined for any positive time; 
\item[(3a)] if $\kkk=\rrr$, there is a smooth function $f:\sss^{n-1}\rightarrow\mathbb R$ such that $\tilde g_t$ converges to the Riemannian metric $\tilde g_{\infty}=e^{2f}\sigma$, moreover there are infinitely many $\mm_0$ such that $\tilde g_{\infty}$ does not have constant scalar curvature;
\item[(3b)] if $\kkk=\ccc$, there is a smooth $\sss^1$-invariant function $f:\sss^{2n-1}\rightarrow\mathbb R$ such that $\tilde g_t$ converges to the sub-Riemannian metric $\tilde g_{\infty}=e^{2f}\sigma_{\ccc}$, moreover there are infinitely many $\mm_0$ such that $\tilde g_{\infty}$ does not have constant Webster scalar curvature;
\item[(3c)] if $\kkk=\hhh$, there is a smooth $\sss^3$-invariant function $f:\sss^{4n-1}\rightarrow\mathbb R$ such that $\tilde g_t$ converges to the sub-Riemannian metric $\tilde g_{\infty}=e^{2f}\sigma_{\hhh}$, moreover there are infinitely many $\mm_0$ such that $\tilde g_{\infty}$ does not have constant quaternionic contact scalar curvature.
\end{itemize}
\end{Theorem}

The main inspiration for this Theorem is a series of papers which considers the evolution of the same class of hypersurfaces, but with respect to the inverse mean curvature flow. The case of the real hyperbolic space has been studied by Gerhard \cite{Ge2}, Ding \cite{Di}, Hung and Wang \cite{HW}. In this last paper the authors showed a fundamental difference with the Euclidean space: in the hyperbolic space there are infinitely many initial data such that the limit is not round. More recently the author of the present paper studied the inverse mean curvature flow in the complex hyperbolic space \cite{Pi4} and in the quaternionic hyperbolic space \cite{Pi3} observing for the first time the presence of a sub-Riemannian limit. In these two papers we showed also that if we wish to classify the possible limits in terms of their curvature, the usual Riemannian scalar curvature is not the right notion because it blows up for any initial datum, therefore we considered the Webster and the quaternionic contact curvature. A survey about the inverse mean curvature flow in the Euclidean and hyperbolic spaces can be read in \cite{Pi2}. 

Since $\psi(x)=x$ satisfies Conditions \eqref{condition}, then Theorem \ref{main} is a generalization to a bigger class of speeds of all this series of results. The presence of the nonhomogeneous speed introduces some technical difficulties. Conditions \eqref{condition} are crucial for proving especially part (1) and (2). In fact i) guarantees the short time existence of the flow, ii) is used for proving that the star-shapeness is preserved and iii) helps for the mean convexity and the long time existence of the flow. In \cite{HW} the authors considered the modified Hawking mass \eqref{def HM} of a hypersurface of $\rrr\hhh^n$. They proved that along the inverse mean curvature flow the limit of the rescaled induced metric is round if and only if this mass converges to zero. Finally they showed that there are infinitely many examples such that their mass does not converges to zero. In Section \ref{conv} we will prove that the modified Hawking mass is an excellent tool for the construction of the examples of Theorem \ref{main} (3a) even for a general $\psi$. If $\kkk\neq\rrr$, the modified Hawking mass cannot be used, therefore in \cite{Pi3,Pi4} we introduced a weaker notion of mass \eqref{def Q} that has the disadvantage that it does not fully classify the initial data with round limit, but it is enough to construct the desired counterexamples. In the final Section we will prove that this mass works very well with the flow \eqref{def flow} too, completing the proof of Theorem \ref{main} (3b) and (3c).

The paper is organized as follows. In Section \ref{pre} we introduce the Webster and the quaternionic contact scalar curvature, we discuss the associated Yamabe problems, we collect some basic notions about the geometry of $\kh^n$ and its hypersurfaces and we list some general properties of our flows, including the evolution equations of the most significant geometric quantities. In Section \ref{1st} we start the proof of Theorem \ref{main} showing that the star-shapeness and the mean convexity are preserved. The main results of Section \ref{2nd} is the long time existence of the flow. The limit of the rescaled induced metric is the topic of Section \ref{conv}: first we will show that, according to the value of $\kkk$, we have the convergence to a conformal multiple of $\sigma_{\kkk}$, then we will pass to the construction of examples which develops a limit with not constant scalar curvature.


\section{Preliminaries}\label{pre}
\subsection{Riemannian and sub-Riemannian metric on the sphere}
Every hypersurface considered in this paper is closed and star-shaped and so it is an embedding of $\sss^m$, the sphere of dimension $m$ into $\rrr^{m+1}\equiv\kkk^n$. On that sphere we will consider different ``standard'' metrics. In the introduction we have already met $\sigma$ (the round metric with constant sectional curvature equal to $1$). When $\kkk=\rrr$ it is the only metric that we need. In the other two cases the richer geometry of $\ccc\hhh^n$ and $\hhh\hhh^n$ imposes to introduce something more. 

If $\kkk\neq\rrr$, we can consider the action of $\sss^a$ on $\sss^m$ that defines the Hopf fibration $\pi:\sss^m\rightarrow\mathbb{KP}^{n-1}$. Let $\mathcal V$ be the \emph{vertical distribution} of $\pi$: it is the distribution tangent to the fibers of $\pi$. Let $\mathcal H$ be its orthogonal (with respect to $\sigma$) complement, it is called the \emph{horizontal distribution}. The canonical deformation of the Hopf fibration produces an important family of Riemannian metrics on $\sss^m$, the so called \emph{Berger metrics}: fix $\kkk\neq\rrr$ and a positive parameter $\lambda$, the Berger metric of parameter $\lambda$ is the metric $e_{\lambda}$
\begin{equation}\label{berger}
e_{\lambda}(X,Y):=\left\{\begin{array}{lcl}
\sigma(X,Y) & \text{if} & X,Y\in\mathcal H;\\
0 & \text{if} & X\in\mathcal H, \ Y\in\mathcal V;\\
\lambda\sigma(X,Y) & \text{if} & X,Y\in\mathcal V;
\end{array}
\right.
\end{equation}
In order not to make the notations too heavy, we use the same symbol ``$e$'' for both the value of $\kkk$, but it is important to keep in mind that \eqref{berger} produces two different family of metrics. In fact, for example, $\dim\mathcal V=a$, hence it depends on $\kkk$. Clearly when $\lambda\rightarrow\infty$ we have that $e_{\lambda}\rightarrow\sigma_{\kkk}$, where $\sigma_{\kkk}$ is the standard sub-Riemannian metric on $\sss^m$ defined in the Introduction. Note that they are sub-Riemannian because in both cases $\mathcal H+[\mathcal H,\mathcal H]$ is the whole tangent space of $\sss^m$, hence $\sigma_{\kkk}$ is enough to define a distance between points on the sphere called the \emph{Carnot-Caratheodory distance}. 

\begin{Notation}\label{notazioni hess}
We introduce the following notation in order to distinguish between derivatives of a function with respect to different metrics. Fix $\kkk\neq\rrr$, for any given function $f:\sss^{m}\rightarrow\mathbb{R}$, let $f_{ij}$ (resp. $\hat{f}_{ij}$ ) be the components of the Hessian of $f$ with respect to $\sigma$ (resp. $e_{\lambda}$). The value of $\lambda$ and $\kkk$ will be clear from the context. The indices go up and down with the associated metric: for instance $\hat{f}_i^k=\hat{f}_{ij}e_{\lambda}^{jk}$, while $f_i^k=f_{ij}\sigma^{jk}$. Analogous notations will be used for higher order derivatives. Moreover here and in the following, unless explicitly stated otherwise, we will always use the Einstein convention about the repeated indices. 
\end{Notation}

For any given function it will be useful to compare its second derivatives computed with respect to the Berger metric and those determined by $\sigma$. This comparison is simpler if we assume the $\sss^a$-invariance. In the following result we summarize Lemma 2.3 of \cite{Pi4} and Lemma 2.3 of \cite{Pi3}.

\begin{Lemma}\label{hessiani}
Fix $\kkk\neq\rrr$, $\lambda>0$ and let $\varphi:\sss^m\rightarrow\rrr$ be a smooth function. If $\varphi$ is $\sss^a$-invariant we have:
\begin{eqnarray*}
\Delta_e\varphi := \hat{\varphi}_i^i & =& \Delta_{\sigma}\varphi:=\varphi_i^i;\\
|\nabla^2_{e}\varphi|^2_{e}:=\hat{\varphi}_i^j\hat{\varphi}^i_j & =& |\nabla^2_{\sigma}\varphi|^2_{\sigma}+2a(\lambda-1)|\nabla_{\sigma}\varphi|^2_{\sigma}=\varphi_i^j\varphi_j^i+2a(\lambda-1)\varphi_i\varphi^i.\\
\end{eqnarray*}
\end{Lemma}

Each one of the metric discussed above carries with it a notion of curvature. When the metric is Riemannian it is obvious what we mean by curvature. On the other hand it can be computed that, as $\lambda\rightarrow\infty$, the sectional curvature of $e_{\lambda}$ diverges. This happens every time we approximate a sub-Riemannian metric with a family of Riemannian metrics. Therefore when we talk about the curvature of $\sigma_{\kkk}$ (and their conformal multiples) we need to clarify what we mean. 

When $\kkk=\ccc$, $\sss^m$ has a in a natural way a $CR$-structure given by the $1$-form $\theta(\cdot)=\sigma(J\nu,\cdot)$, where $\nu$ is unit normal to $\sss^m$ embedded in the standard way in $\rrr^{2n}\equiv\ccc^n$, and $J$ is the complex structure of $\ccc^n$. The sub-Riemannian metric $e^{2f}\sigma_{\ccc}$ can be thought as the restriction to $\mathcal H\times\mathcal H$ of the Webster metric of the $CR$-structure definite by $e^{2f}\theta$. In this context a fundamental notion is the Tanaka-Webster connection which is the unique connection which satisfies some compatibility conditions with the $CR$-structure. Roughly speaking the Levi-Civita connection mixes all the tangent directions, while the Tanaka-Webster ``remembers'' the splitting of the tangent space in the horizontal and vertical distribution. With this connection we can define in the usual formal way a curvature, called \emph{Webster curvature}. It is well known that $\theta$ has constant Webster curvature (equal to $1$), while in general $e^{2f}\theta$ may not. More details and results about $CR$-geometry can be found in the monograph \cite{DT}.

In the same spirit, when $\kkk=\hhh$, $\sss^m$ inherits from $\rrr^{4n}\equiv\hhh^n$ a \emph{quaternionic contact}-structure (qc-structure for short). The role of the Tanaka-Webster connection is played by the Biquard connection. It can be used to define a qc-Ricci tensor and a qc-scalar curvature. Once again the standard qc-structure has constant qc-curvature, but its conformal multiples may not. A good introduction to qc-geometry is \cite{IV}.

A central and classical subject in Geometric Analysis is the Yamabe problem: find, if they exist, the metrics with constant scalar curvature in a fixed conformal class. It has been solved in great generality for all the three concept of curvature mentioned above. The detailed explanation of its solution goes beyond the purposes of the present work. In the next result we focus only on the cases of our interest.

\begin{Lemma}\label{Yamabe}
Let $f:\sss^m\rightarrow\rrr$ be a smooth function. If $\kkk\neq\rrr$ suppose that $f$ is $\sss^a$-invariant. The following characterizations of the solution of the Yamabe problem holds. 
\begin{itemize}
\item[(1)] $e^{2f}\sigma$ has constant scalar curvature if and only if $e^{-f}$ is a linear combination of constants and first eigenfunctions on the sphere;
\item[(2)] $e^{2f}\sigma_{\ccc}$ has constant Webster scalar curvature if and only if $f$ is constant;
\item[(3)] $e^{2f}\sigma_{\hhh}$ has constant qc-scalar curvature if and only if $f$ is constant.
\end{itemize}
\end{Lemma}
Part (1) is the content of Lemma 4 of \cite{HW}, part (2) is Lemma 2.5 of \cite{Pi4} and part (3) is Lemma 2.4 of \cite{Pi3}.

\subsection{Geometry of hyperbolic spaces}
The ambient manifolds that we are considering can be characterized in many ways and they can be described with many different isometric models. Since we wish to work with star-shaped hypersurfaces, the best thing to do is to introduce polar coordinates. It is well known that $\rrr\hhh^n$ can be thought as $\rrr^n$ equipped with the metric
$$
\bar g=d\rho^2+\siq \sigma,
$$
where $\rho$ denotes the radial distance from the center of the coordinates.

In general the underlying manifold of $\kh^n$ is $\rrr^{(a+1)n}\equiv \kkk^n$, where $a$ has been defined in \eqref{def a}, equipped with the metric
$$
\bar g=d\rho^2+\siq e_{\coq}, 
$$
where $e_{\coq}$ is the Berger metric \eqref{berger} of parameter $\coq$. Its curvature tensor has the following explicit expression

\begin{eqnarray}
\nonumber\bar R(X,Y,Z,W) & = &-\bar g(X,Z) \bar g(Y,W)+\bar g(X,W)\bar g(Y,Z)\\
 \label{curv}&&+\sum_{i=1}^a\left[-\bar g(X,J_iZ)\bar g(Y,J_iW)+\bar g(X,J_iW)\bar g(Y,J_iZ)\right]\\
\nonumber&&-2\sum_{i=1}^a\bar g(X,J_iY)\bar g(Z,J_iW),
\end{eqnarray}
where $J_1,\dots,J_a$ are the complex structure of $\kh^n$ that, in this model, coincide with those of $\kkk^n$. Note that if $\mathbb K=\mathbb R$, and hence $a=0$, the sums in the second and third line of \eqref{curv} are empty. Therefore \eqref{curv} can be used to describe the curvature tensor of $\rrr\hhh^n$ too.

From \eqref{curv} it follows that the our ambient manifolds are symmetric, the curvature is constant equal to $-1$ if $\kkk=\rrr$, otherwise is bounded between $-4$ and $-1$. Moreover $\kh^n$ is Einstein with Ricci tensor given by

\begin{equation}\label{ricci}
\bar Ric =-(m+3a)\bar g.
\end{equation}

\subsection{Geometry of hypersurfaces in hyperbolic spaces}
Let $\mm$ be a real closed star-shaped hypersurface of $\kh^n$, then it is an embedding of the sphere of dimension $m$ in $\kh^n$. Up to an isometry of the ambient manifold we can always suppose that it is star-shaped with respect to the center of the polar coordinates.  Hence there is a positive function $\rho:\sss^m\rightarrow\mathbb R$ such that in polar coordinate $\mm=\left\{(x,\rho(x))\in\kh^n\ \left|\ x\in\sss^m\right.\right\}$. Vice versa each positive function defines a star-shaped hypersurface via the embedding
$$
\begin{array}{cccc}
F:&\sss^m&\rightarrow&\kh^n\\
& x& \mapsto& (x,\rho(x)).
\end{array}
$$
For any given $\mm$, we call such $\rho$ the \emph{radial function} associated to $\mm$. If $\mathbb K\neq\mathbb R$, we know that $\sss^a$ acts by isometries on $\sss^m$. In this case we say that $\mm$ is $\sss^a$\emph{-invariant} if its radial function is $\sss^a$-invariant. For reasons of synthesis, sometimes will talk about, with an abuse of notation, $\sss^a$-invariance even when $\kkk=\rrr$: in this case it has to be considered as an empty condition. With the same proof of Lemma 3.1 of \cite{Pi1} we can prove the following result.

\begin{Lemma} The evolution of an $\sss^a$-invariant hypersurface of $\kh^n$ stays $\sss^a$-invariant during the whole duration of the flow.
\end{Lemma}

Now we want to describe the main geometric quantities associated to a star-shaped hypersurface in term of its radial function. For more details and explicit computations we refer to Section 3 of \cite{Pi3} and Section 3 of \cite{Pi4}.

For technical reasons we introduce an auxiliary function $\varphi=\varphi(\rho)$ such that $\frac{d\varphi}{d\rho}=\frac{1}{\sih}$. Fix $(Y_1,\dots,Y_m)$ a tangent basis of $\sss^m$ and denote with $\rho_i:=Y_i(\rho)$ and with $\varphi_i:=Y_i(\varphi)=\frac{\rho_i}{\sih}$. When $\mathbb K\neq\mathbb R$, since the ambient metric is no more isotropic, it is convenient to choose a tangent basis on $\sss^m$ adapted to the contact structure: from now on we always suppose that for any $i=1,\dots,a$ $Y_i=J_i\nu$, where $\nu$ is the unit normal of the standard immersion of $\sss^m$ in $\mathbb R^{m+1}\equiv \mathbb K^n$. The use of this base simplifies some computations because, for example, $\rho_i=\varphi_i=0$ for every $i=1,\dots,a$. Let $V_i:=F_*Y_i=Y_i+\rho_i\partial\rho$, then $(V_1,\dots,V_m)$ is a basis of the tangent space of $\mm$. Let $g=F^*\bar g$ be the induced metric on $\mm$. In coordinates it can be expressed as
\begin{equation}\label{ind metr}
g_{ij}=\siq\left(\varphi_i\varphi_j+e_{ij}\right).
\end{equation}

\begin{Notation}\label{nnn}
For the proofs in the next Sections it will be useful to have a common notation that allows us to work with all values of $\kkk$ at once. Therefore here and in the sequel often we will use the symbol ``$e$'' even if $\kkk=\rrr$. In this case the Berger metrics are not defined and $e$ is just another way to denote $\sigma$. For example, with this notation, \eqref{ind metr} holds for any value of $\kkk$. In the same spirit when $\kkk=\rrr$ the symbol ``$\sigma_{\kkk}$'' denotes again $\sigma$. For clarity we specify that on the other hand the symbol ``$\sigma$'' is reserved uniquely for the round Riemannian metric on the sphere. 
\end{Notation}

\noindent The outward unit normal vector field of $\mm$ is
\begin{equation}\label{normal}
\nu=\frac{1}{v}\left(\partial\rho-\frac{\nabla\varphi}{\sih}\right),
\end{equation}
where
$$
v=\bar g(\nu,\partial\rho)^{-1}=\sqrt{1+|\nabla\varphi|^2}.
$$
The gradient $\nabla$ should be taken with respect to $\sigma$ if $\mathbb K=\mathbb R$ and with the Berger metric otherwise, but since we are considering only $\sss^a$-invariant hypersurfaces we can consider $\sigma$ in any case.
The inverse of the induced metric is
$$
g^{ij}=\frac{1}{\siq}\left(e^{ij}-\frac{\varphi^{i}\varphi^{j}}{v^2}\right),
$$
where we are using Notation \ref{nnn} when $\kkk=\rrr$. The second fundamental form of $\mm$ is

\begin{equation}\label{2ff}
h_i^j=-\frac{\hat{\varphi}_{ik}\tilde e^{kj}}{v\sih}+\frac{\coh}{v\sih}\delta_i^j+\frac{\sih}{v\coh}\sum_{k=1}^a\delta_i^k\delta_k^j,
\end{equation}

\noindent where $\tilde e^{ij}=\siq g^{ij}=e^{ij}-\frac{\varphi^i\varphi^j}{v^2}$ and we used Notation \ref{notazioni hess} for the second derivative of $\varphi$. Taking the trace of \eqref{2ff} and using Lemma \ref{hessiani} we can compute the mean curvature of $\mm$:
\begin{equation}\label{H}
H=-\frac{\varphi_{ij}\tilde{\sigma}^{ij}}{v\sih}+\frac{\hat H}{v},
\end{equation}
where $\tilde{\sigma}^{ij}=\sigma^{ij}-\frac{\varphi^i\varphi^j}{v^2}$ and 
\begin{equation}\label{H hat}
\hat H(\rho)=m\cth+a\tgh.
\end{equation}

\subsection{Evolution equations}

Since by Condition \eqref{condition} we have that $\psi'>0$, the short-time existence and uniqueness of the solution for the flow \eqref{def flow} are guaranteed by standard arguments. Moreover well known computation (see for example \cite{HP}) can be repeated to compute the evolution equation of the main geometric quantities.

\begin{Lemma}\label{eq evoluz} Since the ambient space is symmetric the following evolution equations hold:
\begin{itemize}
\item[(1)] $\displaystyle{\frac{\partial g_{ij}}{\partial t} = \frac {2}{\psi} h_{ij}},\qquad \displaystyle{\frac{\partial g^{ij}}{\partial t} = -\frac {2}{\psi} h^{ij}}$,
\item[(2)] $\displaystyle{\frac{\partial H}{\partial t} = \frac{\psi'}{\psi^2}\Delta H+\frac{\psi''\psi-2(\psi')^2}{\psi^3}\left|\nabla H\right|^2-\frac{1}{\psi}\left(\left|A\right|^2+\bar{R}ic(\nu,\nu)\right)}$,
\item[(3)] $\displaystyle{\frac{\partial \psi}{\partial t} = \frac{\psi'}{\psi^2}\Delta \psi-2\frac{(\psi')^3}{\psi^3}\left|\nabla H\right|^2-\frac{\psi'}{\psi}\left(\left|A\right|^2+\bar{R}ic(\nu,\nu)\right)}$,
\item[(4)] $\displaystyle{\frac{\partial h_i^j}{\partial t} =-\nabla_i\nabla^j\frac{1}{\psi}+\frac{1}{\psi}\left(\bar R_{0i0}^{\phantom{0i0}j}-h_i^kh_k^j\right)}$,
\item[\phantom{(4)}]$\displaystyle{\phantom{\frac{\partial h_i^j}{\partial t}}=\frac{\psi'}{\psi^2}\Delta h_i^j+\frac{\psi''\psi-2(\psi')^2}{\psi^3}\nabla_iH\nabla^jH-\left(\frac{1}{\psi}+\frac{H\psi'}{\psi^2}\right)\left(h_i^kh_k^j+\bar R_{0i0}^{\phantom{0i0}j}\right)}$
\item[\phantom{(4)}]$\displaystyle{\phantom{\frac{\partial h_i^j}{\partial t} =}+\frac{\psi'}{\psi^2}\left(\left(\left|A\right|^2+\bar{R}ic(\nu,\nu)\right)h_i^j+2\bar R_{\phantom{k}is\phantom{j}}^{k\phantom{is}j}h_k^s-\bar R_{\phantom{k}si\phantom{s}}^{k\phantom{si}s}h_k^j-\bar R_{\phantom{s}ks\phantom{j}}^{s\phantom{ks}j}h_i^k\right)}$,
\item[(5)] $\displaystyle{\frac{\partial d\mu}{\partial t}=\frac{H}{\psi}d\mu}$.
\end{itemize}
\end{Lemma}
Here and in the following, if there is no risk of confusion, we are using for brevity only $\psi$ for saying $\psi(H)$, and analogously for its derivatives. Note that integrating the evolution equation of the volume form $d\mu$ we get
\begin{equation}\label{ev vol}
\frac{d\left|\mm_t\right|}{d t}= \int_{\mm_t}\frac{H}{\psi}d\mu, 
\end{equation}
in particular it follows that \eqref{def flow} is an expanding flow, at least as far the evolving hypersurface is mean convex.

\begin{Example}
A geodesic sphere is a star-shaped hypersurface with constant radial function. Therefore, by \eqref{H} its mean curvature is given by $\hat H(\rho)$. In particular it is constant. Therefore it is easy to see that the evolution of a geodesic sphere is a family of geodesic spheres such that the radius evolves in the following way
$$
\frac{d \rho}{d t} = \frac{1}{\psi(\hat H(\rho))}.
$$
For a general $\psi$ we cannot find the explicit solution of this ODE, but $\frac{\partial\hat H}{\partial \rho}<0$ by direct computations, $\psi'>0$ by Conditions \eqref{condition}, therefore $\rho$ is increasing and it blows up in infinite time. Since $\lim_{\rho\rightarrow +\infty}\hat H(\rho)=m+a$, then 
$$
\rho\approx \frac{t}{\psi(m+a)},\ \text{as}\ t\rightarrow +\infty.
$$
\end{Example}


\section{First order estimetes}\label{1st}
The main goal of this Section is to prove that the initial conditions are preserved, i.e. part (1) of Theorem \ref{main}. The most important technical result is the following.

\begin{Proposition}\label{grad exp}
There is a positive constant $c$ such that if $\psi$ satisfies Conditions \eqref{condition} i) and ii), then 
\begin{itemize}
\item[(1)] $\left|\nabla\varphi\right|^2\leq ce^{-\frac{2}{\pma}t}$;
\item[(2)] $|H-m-a|\leq ce^{-\frac{2}{\psi(m+a)}t},\quad|\psi(H)-\psi(m+a)|\leq ce^{-\frac{2}{\pma}t}.$
\end{itemize}
\end{Proposition}

\noindent Part (1) of this Proposition has a direct important consequences.

\begin{Corollary} \label{stellato}
The evolution o
f any star-shaped $\sss^a$-invariant hypersurface stays star-shaped for any time the flow is defined.
\end{Corollary}
\proof By Proposition \ref{grad exp} there exists a positive constant $c$ such that
$$
v=\bar g\left(\frac{\partial}{\partial\rho},\nu\right)^{-1}=\sqrt{1+|\nabla\varphi|^2}\leq c.
$$
It follows that $\frac{\partial}{\partial\rho}$ and $\nu$ are never orthogonal in $\kh^n$. This means that $\mm_t$ is star-shaped for any time $t$.
\cvd

The proof of Proposition \ref{grad exp} proceeds by steps: first we prove that $\nabla\varphi$ is just bounded (which is already enough for having Corollary \ref{stellato}), then we prove that it decays exponentially fast, finally we find the optimal exponent. In the meanwhile we are able to show that $H$ stays strictly positive and bounded, and converges exponentially fast to $m+a$, i.e. to the mean curvature of a horosphere in $\kh^n$, finally we can find the optimal exponent for $H$ too. The first crucial step is the following Lemma.

\begin{Lemma}\label{grad lim}
If $\psi$ satisfies Conditions \eqref{condition} i) and ii), then for any $(x,t)$
$$
\left|\nabla\varphi(x,t)\right|^2\leq \max_{y\in\sss^m}\left|\nabla\varphi(y,0)\right|^2.
$$
\end{Lemma}
\proof
Let us define $\omega=\frac 12|\nabla\varphi|^2=\frac 12\varphi_k\varphi^k$. We want to compute the evolution equation of $\omega$ and apply the maximum principle. We have that the radial function satisfies the scalar evolution equation $\frac{\partial\rho}{\partial t} =\frac{v}{\psi(H)}$, hence the evolution of $\varphi$ is given by 
\begin{equation}\label{ev phi}
\frac{\partial\varphi}{\partial t} = G(\varphi_{ij},\varphi_i,\varphi):=\frac{v}{\sih \psi(H)}.
\end{equation}
The original geometric flow \eqref{def flow} is defined at least as far the scalar flow \eqref{ev phi} is defined and, when both are defined, they are equivalent. Therefore we can work with \eqref{ev phi}.

Let $a^{ij}=\frac{\partial G}{\partial \varphi_{ij}}  =  \frac{\psi'}{\psi^2} g^{ij}$: it is a symmetric and positive definite, at least as far $\frac{\psi'}{\psi^2}$ is bounded and strictly positive. Moreover we denote by $b^i=\frac{\partial G}{\partial \varphi_{i}}$. From \eqref{ev phi} we have:
\begin{eqnarray*}
\frac{\partial\omega}{\partial t} & = & \varphi^k\nabla_k\frac{\partial\varphi}{\partial t}=\varphi^k\left(a^{ij}\varphi_{ijk}+b^i\varphi_{ik}+\frac{\partial G}{\partial\varphi}\varphi_k\right).
\end{eqnarray*}
Let $R$ be the Riemannian curvature tensor of $\sigma$, then the Ricci identity says
\begin{equation}\label{r id}
\varphi_{ijk}=\varphi_{kij}+R^{m}_{\phantom{m}ijk}\varphi_m=\varphi_{kij}+\varphi_j\sigma_{ik}-\varphi_k\sigma_{ij}.
\end{equation}
Since $a^{ij}$ is symmetric and positive definite, applying \eqref{r id}, after some explicit computations we get
$$
a^{ij}\varphi_{ijk}\varphi^k = a^{ij}\omega_{ij}-a^{ij}\varphi_{ik}\varphi_j^k+a^{ij}\varphi_i\varphi_j-2a^i_i\omega \leq a^{ij}\omega_{ij}.
$$
As a consequence of the $\sss^a$-invariance we have
$$
\frac{\partial\tilde{\sigma}}{\partial\varphi}=0,\quad\frac{\partial v}{\partial\varphi}=0,
$$
hence, using the explicit expression of the mean curvature \eqref{H} we can compute
\begin{eqnarray}
\label{G_phi precisa}\frac{\partial G}{\partial\varphi}& = & \frac{v\coh}{\sih\psi}\left(\frac{H\psi'}{\psi}-1\right)-\frac{\psi'}{\psi^2}\left(m+a+\frac{a}{\coq}\right)\\
& \leq &-\frac{\psi'}{\psi^2}\left(m+a+\frac{a}{\coq}\right) \label{G_phi}
\end{eqnarray}
where in the last line we used Condition \eqref{condition} ii). Summarizing we have just found that
$$
\frac{\partial\omega}{\partial t}\leq a^{ij}\omega_{ij}+b^i\omega_i.
$$
The result follows by the maximum principle.
\cvd

Now we are able to prove that $H$ is strictly positive and bounded.

\begin{Lemma}\label{H lim}
If $\psi$ satisfies Conditions \eqref{condition}, then there exist two positive constant $c_1,c_2$ such that for any time the flow is defined we have
$$
0<c_1\leq H\leq c_2.
$$
\end{Lemma}
\proof
We can start from the upper bound. Combining the evolution equation of $H$ given in Lemma \ref{eq evoluz}, with Condition \eqref{condition} iii), the fact that $H^2\leq m |A|^2$ and \eqref{ricci} we have:
\begin{eqnarray*}
\frac{\partial H}{\partial t}&\leq &\frac{\psi'}{\psi^2}\Delta H-\frac{1}{\psi}\left(\frac{H^2}{m}-m-3a\right).
\end{eqnarray*}
By the maximum principle we have that $H\leq c_2$ for some constant $c_2$ depending only on $m$, $a$ and $\mm_0$. On the other hand, let $\alpha=\frac{\partial\varphi}{\partial t}=\frac{v}{\sih\psi}$, then by \eqref{G_phi} and \eqref{ev phi}
\begin{eqnarray*}
\frac{\partial \alpha}{\partial t} & = & a^{ij}\alpha_{ij}+b^i\alpha_i+\frac{\partial G}{\partial\varphi}\alpha\leq a^{ij}\alpha_{ij}+b^i\alpha_i.
\end{eqnarray*}
By the maximum principle $\alpha$ is bounded from above, therefore there is a positive constant $c$ such that
$$
\psi\geq \frac{c}{\sih}.
$$
Since $\rho$ does not blow up in finite time, this means that $\psi$, and hence $H$, are strictly positive for any finite time. Now we can improve what just said showing that $H$ cannot converge to zero. Let us consider the function $\tilde{\alpha}=\frac{v}{\sih\psi}e^{\frac {t}{\pma}}$. Recalling that $\rho\approx\frac {t}{\pma}$ as $t\rightarrow +\infty$ ( if we can take arbitrary big times $t$), Lemma \ref{grad lim}, and the fact that $\psi'>0$, then an upper bound for $\tilde{\alpha}$ implies a strictly positive lower bound for $H$. We compute the evolution equation of $\tilde{\alpha}$:
\begin{eqnarray}
\frac{\partial \tilde{\alpha}}{\partial t} & = & a^{ij}\tilde{\alpha}_{ij}+b^i\tilde{\alpha}_i+\frac{\partial G}{\partial\varphi}\tilde{\alpha}+\frac {1}{\pma}\tilde{\alpha}\label{ev tilde alpha}\\
\nonumber &\leq&a^{ij}\tilde{\alpha}_{ij}+b^i\tilde{\alpha}_i-\tilde{\alpha}(m+a)\frac{\psi'}{\psi^2}+\frac {1}{\pma}\tilde{\alpha}\\
\nonumber & = & a^{ij}\tilde{\alpha}_{ij}+b^i\tilde{\alpha}_i-\tilde{\alpha}^2(m+a)\frac{\psi'\sih}{ve^{\frac{t}{\pma}}\psi}+\frac {1}{\pma}\tilde{\alpha}.
\end{eqnarray}
Since $\psi(0)=0$ then there exists a positive constant $c$ such that $\frac{\psi'(H)}{\psi(H)}\geq c$ for every $H\in[0,c_2]$ even if $\lim_{x\rightarrow 0}\psi'(x)=0$. Therefore 
\begin{eqnarray*}
\frac{\partial \tilde{\alpha}}{\partial t} & \leq & a^{ij}\tilde{\alpha}_{ij}+b^i\tilde{\alpha}_i-\tilde c\tilde{\alpha}^2+\frac {1}{\pma}\tilde{\alpha},
\end{eqnarray*}
for some $\tilde c>0$. By the maximum principle we can conclude that $\tilde{\alpha}$ is bounded from above from a constant that does not depend on time.
\cvd

Now we can improve what said so far showing that $|\nabla\varphi|^2$ decays exponentially fast and that $H$ converges (once we will prove long time existence for the flow) exponentially fast to $m+a$.

\begin{Lemma}\label{H exp}
There exist positive constants $\beta,\ \gamma,\ c$ such that 
\begin{itemize}
\item[(1)] $|\nabla\varphi|^2\leq ce^{-\beta t}$;
\item[(2)] $|H-m-a|\leq ce^{-\gamma t};\quad |\psi(H)-\psi(m+a)|\leq ce^{-\gamma t}$.
\end{itemize}
\end{Lemma}
\proof
\begin{itemize}
\item[(1)] The function $\frac{\psi'}{\psi^2}$ is continuous and strictly positive in $[c_1,c_2]$, where $c_1$ and $c_2$ are determined in Lemma \ref{H lim}. Then we can find a strictly positive constant $b$ such that for any time $t$
$$
\frac{\psi'}{\psi^2}\geq b.
$$
Hence we can define $\beta=(m+a)b$ and repeat the proof of Lemma \ref{grad lim} improving the estimates of the reaction term in order to have an exponential decay:
$$
\frac{\partial G}{\partial \varphi}\leq-(m+a)\frac{\psi'}{\psi^2}\leq -\beta.
$$
\item[(2)] When $\mathbb K=\mathbb R$ we have $a=0$ and $H$ can be estimate as follows:
\begin{eqnarray*}
\frac{\partial H}{\partial t} &\leq & \frac{\psi'}{\psi^2}\Delta H-\frac{1}{\psi}\left(|A|^2-m\right)\leq \frac{\psi'}{\psi^2}\Delta H-\frac{1}{m\psi}\left(H^2-m^2\right).
\end{eqnarray*}
Since in Lemma \ref{H lim} we proved that $H$ is bounded and it cannot be too close to $0$, we have that $\frac{1}{\psi}$ is bounded too. Therefore applying the maximum principle we get that there exists a constant $\gamma>0$ such that
$$
H-m\leq c e^{-\gamma t}
$$
When $\mathbb K\in\left\{\mathbb C,\mathbb H\right\}$, the proof is more involved. By the $\sss^a$-invariance, for any $k=1,\dots,a$ \eqref{2ff} reduced to $h_k^k=\frac{\coh}{v\sih}+\frac{\sih}{v\coh}$. By part (1) of this Lemma we have that for any $k=1,\dots,a$
$$
|h_k^k-2|\leq ce^{-\beta t}.
$$
Arguing as in Lemma 5.6 of \cite{Pi4}, we can define the tensor $l_i^j:=h_i^j+\delta_i^j-\sum_{k=1}^a \delta_i^k\delta_k^j,$ and its trace $L:=l_i^i=H+m-a$. We get:
\begin{eqnarray*}
\frac{\partial H}{\partial t} & \leq & \LLL H-\frac{1}{\psi}\left(|A|^2+\bar Ric(\nu,\nu)\right)\\
& = & \LLL H-\frac{1}{\psi}\left(|l|^2-2(H-m-a)+2\sum_{k=1}^a(h_k^k-2)-4m\right)\\
& \leq & \LLL H-\frac{1}{\psi}\left(\frac{L^2}{m}-2(H-m-a)-4m\right)+ce^{-\beta t}\\
& = & \LLL H-\frac{1}{m\psi}\left(H-m-a\right)(H+m-a)+ce^{-\beta t}.
\end{eqnarray*}
Since $\frac{H+m-a}{m\psi}$ is strictly positive and, by Lemma \ref{H lim}, bounded, by the maximum principle there is a constant $0<\gamma\leq\beta$ such that
\begin{equation}\label{H1}
H-m-a\leq ce^{-\gamma t}.
\end{equation}

On the other hand, the strategy for proving the lower bound is similar to that of Lemma 5.6 of \cite{Pi4}, however we need some modifications in order to consider a general $\psi$. In fact, like in the proof of Lemma \ref{H lim}, we are able to estimate $H$ from below only considering first $\psi$ with the help of the function $\tilde{\alpha}=\frac{ve^{\frac{t}{\pma}}}{\psi\sih}$. We restart from \eqref{ev tilde alpha}, the evolution equation, the n of $\tilde{\alpha}$. This time, we need to use a finer estimate on the reaction term: from \eqref{G_phi precisa}, part (1) of this Lemma and \eqref{H1} we get
\begin{eqnarray*}
\frac{\partial \tilde{\alpha}}{\partial t} & \leq & a^{ij}\tilde{\alpha}_{ij}+b^i\tilde{\alpha}_i+\frac{\tilde{\alpha}}{\pma}+\left(\frac{\psi'}{\psi^2}\left(\frac{v\coh}{\sih}H-m-a\right)-\frac{v\coh}{\psi\sih}\right)\tilde{\alpha}\\
& \leq &a^{ij}\tilde{\alpha}_{ij}+b^i\tilde{\alpha}_i+\left(\frac {1}{\pma}+ce^{-\gamma t}-\frac{\tilde{\alpha}}{2}\right)\tilde{\alpha}.
\end{eqnarray*}
Therefore by the maximum principle $\tilde{\alpha}\leq\frac{2}{\pma}+ce^{-\gamma t}$, hence, by definition of $\tilde{\alpha}$ we have
\begin{equation}\label{H2}
\psi(H)-\psi(m+a)\geq -ce^{-\gamma t}.
\end{equation}
We can combine \eqref{H1} and \eqref{H2} to get the result using the mean value theorem.
\end{itemize}
\cvd

\noindent Finally we can look for the optimal exponent.

\noindent\emph{Proof of Proposition \ref{grad exp}.}
\begin{itemize}
\item[(1)] As in Lemma \ref{grad lim} we consider again the evolution of $w$. Using Lemma \ref{H exp} and estimating the reaction terms of the evolution equation of $\omega$ with \eqref{G_phi precisa} we get
$$
\frac{\partial w}{\partial t} \leq a^{ij}w_{ij}+b^iw_i+\left(ce^{-\gamma t}-\frac{2}{\pma}\right)w.
$$
The result follows from the maximum principle.
\item[(2)] The proofs of Lemma \ref{H exp}, part (2) can be repeated using $\beta=\frac{2}{\pma}$ and noting that $\frac{H+m-a}{m\psi}$ converges exponentially fast to $\frac{2m-a}{m\pma}$ which is smaller than $\frac{2}{\pma}$, therefore we can take $\lambda=\frac{2}{\pma}$.
\end{itemize}
\cvd


\section{Higher order estimate and long time existence}\label{2nd}
The main goal of this section is to prove part (2) of Theorem \ref{main}, i.e. the long time existence of the flow. Moreover we will show some other important auxiliary results, such as the convergence of the second fundamental form to that of a horosphere of $\kh^n$.

\begin{Proposition}\label{EEE}
The principal curvatures of the evolving hypersurface are uniformly bounded for any time.
\end{Proposition}
\proof
Since by the results of the previous section $H$ is bounded from below, it is sufficient to prove that the principal curvatures are bounded from above. We want to adapt the strategy of Corollary 6.2 of \cite{Pi4} to the case of a nonhomogeneous flow, therefore we define the tensor $M_i^j=\psi(H) h_i^j$. By Lemma \ref{eq evoluz} and after some standard computations we have that the evolution equation of $M_i^j$ is
\begin{eqnarray*}
\frac{\partial M_i^j}{\partial t} & = & \frac{\psi'}{\psi^2}\Delta M_i^j-2\frac{\psi'}{\psi^3}\left\langle \nabla M_i^j,\nabla\psi\right\rangle+\frac{\psi''\psi-2(\psi')^2}{\psi^2}\nabla_iH\nabla^jH\\
& &-\frac{1}{\psi^2}\left(1+\frac{H\psi'}{\psi}\right)M_i^kM_k^j-\left(1+\frac{H\psi'}{\psi}\right)\bar R_{0i0}^{\phantom{0i0}j}\\
&&+\frac{\psi'}{\psi^2}\left(2\bar R_{\phantom{k}is\phantom{j}}^{k\phantom{is}j}M_k^s-\bar R_{\phantom{k}si\phantom{s}}^{k\phantom{si}s}M_k^j-\bar R_{\phantom{s}ks\phantom{j}}^{s\phantom{ks}j}M_i^k\right).
\end{eqnarray*}
Let $\mu_1\leq\mu_2\leq\dots\leq\mu_m$ be the eigenvalues of $M_i^j$. Since $\sum_{i=1}^m\mu_i=H\psi(H)>0$, then $\mu_m>0$ everywhere. The goal is to prove that $\mu_m$ is bounded from above. Fix any time $T^*$ strictly smaller than the maximal time $T$ . We can find a point $(x_0,t_0)$ where $\mu_m$ reaches its maximum in $\sss^m\times [0,T^*]$. At this point we can fix an orthonormal basis
which diagonalizes $M_i^j$, then we can say that at this point $\mu_m$ satisfies the same evolution equation of $M_m^m$. By Condition \eqref{condition} the term $\frac{\psi''\psi-2(\psi')^2}{\psi^2}\nabla_mH\nabla^mH$ is negative and it can be ignored. Moreover the curvature of $\kh^n$ is bounded, hence all the coefficients involving $\bar R$ are bounded. Furthermore $H,\ \psi$ and $\psi'$ are uniformly bounded too, therefore there are positive constants $C_0,\ C_1,\ C_2$ independent on the choice of $T^*$ such that in $(x_0,t_0)$ the following holds
$$
\frac{\partial\mu_m}{\partial t}\leq \frac{\psi'}{\psi^2}\Delta \mu_m-2\frac{\psi'}{\psi^3}\left\langle \nabla \mu_m,\nabla\psi\right\rangle-C_2\mu_m^2+C_1\mu_m+C_0.
$$
Since $\mu_m$ is positive, it follows that it is bounded by a constant depending on $\mm_0$, but not on the choice of $T^*$.
\cvd

\begin{Corollary}\label{volume}
The flow is definite for any positive time. Moreover the flow is expanding and, as $t$ diverges, the volume goes like
$$
|\mm_t|\approx e^{\frac{m+a}{\psi(m+a)}t}.
$$
\end{Corollary}
\proof
The long time existence of the flow is an immediate consequence of Proposition \ref{EEE}. For the growth of the volume, consider the quantity $V=|\mm_t|e^{-\frac{m+a}{\psi(m+a)}t}$. By Lemma \ref{eq evoluz} and Proposition \ref{grad exp} we have
\begin{eqnarray*}
\left|\frac{d V}{dt}\right|& = & e^{-\frac{m+a}{\psi(m+a)}t}\left|\int_{\mm_t}\left(\frac{H}{\psi}-\frac{m+a}{\psi(m+a)}\right)d\mu_t\right|\leq ce^{-\frac{2}{\psi(m+a)}t}V.
\end{eqnarray*}
Hence $0<|\mm_0|e^{-\frac{c\pma}{2}}\leq V\leq |\mm_0|e^{\frac{c\pma}{2}}$.
\cvd


\begin{Lemma}\label{higher}
For every $N\in\mathbb N$ there is a positive constant $c$ such that the $N$-th derivative of $\varphi$ satisfies
$$
|\nabla^N\varphi|^2\leq ce^{-\frac{2}{\pma}t}.
$$
\end{Lemma}
\proof
Fix $N\in\mathbb N$ and consider the quantity $\omega=\frac{1}{2}|\nabla^N\varphi|^2=\frac 12\varphi^{k_1\dots k_N}\varphi_{k_1\dots k_N}$.
With the notations of the proof of Lemma \ref{grad lim}, its evolution equation is
\begin{eqnarray*}
\frac{\partial\omega}{\partial t} & = & \varphi^{k_1\dots k_N}\nabla_{k_1}\dots\nabla_{k_N}\frac{\partial\varphi}{\partial t}\\
&=&\varphi^{k_1\dots k_N}\left(a^{ij}\varphi_{ijk_1\dots k_N}+b^i\varphi_{ik_1\dots k_N}\right)+2\frac{\partial G}{\partial \varphi}\omega
\end{eqnarray*}
By \eqref{G_phi precisa} and Proposition \ref{grad exp} we have that $\frac{\partial G}{\partial\varphi}\leq C_1 e^{-\frac{2}{\pma}t}-\frac{1}{\pma}$ for some constant $C_1>0$ . Applying a finite number of times the Ricci identity \eqref{r id} we have
\begin{eqnarray*}
\frac{\partial\omega}{\partial t} & \leq & a^{ij}\omega_{ij}+b^i\omega_i+\left(2C_1e^{-\frac{2}{\pma}t}-\frac{2}{\pma}\right)\omega\\
&&+a\ast\nabla^N\varphi\ast\nabla^N\varphi+b\ast\nabla^{N-1}\varphi\ast\nabla^N\varphi,
\end{eqnarray*}
where given two tensors $S$ and $T$, $S\ast T$ denotes any linear combination obtained contracting $S$ and $T$ by $\sigma$. We have that
\begin{equation}\label{02}
a^{ij}=\frac{\psi'}{\psi^2}g^{ij}=\frac{\psi'}{\siq\psi^2}\tilde{\sigma}^{ij}\leq C_2e^{-\frac{2}{\pma}t}\sigma^{ij}.
\end{equation}
Moreover, by direct computations we have that
\begin{equation}\label{03}
b^i=\frac{1}{\sih}\nabla\varphi\ast\nabla^2\varphi.
\end{equation}
The case $N=1$ has been already proved in Proposition \ref{grad exp}. Suppose by induction that the result holds for $N-1$, then by \eqref{02} and \eqref{03} we have that there is a positive constant $C_3$ such that
$$
a\ast\nabla^N\varphi\ast\nabla^N\varphi+b\ast\nabla^{N-1}\varphi\ast\nabla^N\varphi\leq C_3e^{-\frac{2}{\pma}t}\omega.
$$
Therefore there is a $C_4>0$ such that
$$
\frac{\partial\omega}{\partial t} \leq  a^{ij}\omega_{ij}+b^i\omega_i+\left(C_4 e^{-\frac{2}{\pma}t}-\frac{2}{\pma}\right)\omega.
$$
The desired estimates follow by the maximum principle.
\cvd

A consequence of this Lemma is the convergence of the second fundamental form to that of a horosphere.

\begin{Corollary}\label{conv 2ff}
There is a positive constant $c$ such that
\begin{itemize}
\item[(1)] if $\kkk=\rrr$ we have
\begin{eqnarray*}
\left|h_i^j-\delta_i^j\right|^2\leq c e^{-\frac{4}{\pmm}t},&&\aao\leq c e^{-\frac{4}{\pmm}t};
\end{eqnarray*}
\item[(2)] if $\kkk\neq\rrr$ we have 
$$
\left|h_i^j-\delta_i^j-\sum_{k=1}^a\delta_i^k\delta_k^j\right|^2\leq c e^{-\frac{2}{\pma}t},
$$
while on the horizontal distribution we have a faster convergence
$$
\sum_{i,j=a+1}^m(h_i^j-\delta_i^j)(h_j^i-\delta_j^i)\leq c e^{-\frac{4}{\pma}t}.
$$
\end{itemize}
\end{Corollary}
\proof
\begin{itemize}
\item[(1)] From the explicit expression of the second fundamental form \eqref{2ff} we get
\begin{eqnarray*}
\left|h_i^j-\delta_i^j\right|^2&=& \frac{1}{v^2\siq}\left|\varphi_{ik}\tilde{\sigma}^{kj}\right|^2+m\left(\frac{\coh}{v\sih}-1\right)^2\\
&&+2\left(\frac{\coh}{v\sih}-1\right)\left(H-\frac{\hat H}{v}\right)\leq ce^{-\frac{4}{\pmm}t},
\end{eqnarray*}
where, in the inequality, we used Proposition \ref{grad exp} and Lemma \ref{higher} with $N=2$. Moreover by the inequality that we have just proved and Proposition \ref{grad exp} 
$$
\aao=\left|h_i^j-\frac Hm\delta_i^j\right|^2\leq \left|h_i^j-\delta_i^j\right|^2+m\left(\frac Hm-1\right)^2\leq ce^{-\frac{4}{\pmm}t}.
$$
\item[(2)] By \eqref{2ff}, Lemma \ref{hessiani}, Proposition \ref{grad exp}, Proposition \ref{grad exp} and Lemma \ref{higher} we have
\begin{eqnarray*}
\left|h_i^j-\delta_i^j-\sum_{k=1}^a\delta_i^k\delta_k^j\right|^2&=& \frac{1}{v^2\siq}\left|\varphi_{ik}\tilde{\sigma}^{kj}\right|^2\\
&&+\frac{2a}{v^2}|\nabla\varphi|^2+2\left(\frac{\coh}{v\sih}-1\right)\left(H-\frac{\hat H}{v}\right)\\
&&+m\left(\frac{\coh}{v\sih}-1\right)^2+a\left(\frac{\sih}{v\coh}-1\right)^2\\
&&+2a\left(\frac{\coh}{v\sih}-1\right)\left(\frac{\sih}{v\coh}-1\right)\\
&\leq&ce^{-\frac{2}{\pma}t}.
\end{eqnarray*}
On the horizontal distribution the computations are similar to those of part (1) of this Corollary, hence we have a faster convergence.
\end{itemize}
\cvd


\section{Convergence and curvature of the induced metric}\label{conv}

We finish this paper with the proof of statements (3a), (3b) and (3c) of Theorem \ref{main}. This is the part where the geometries of the ambient manifolds influence mostly the result. Each geometry produces its own typical behaviour, and all of them are very different from what found in the Euclidean case in \cite{CGT}. We recall that in the Euclidean case the limit is always $\sigma$. We will show in a while that in the hyperbolic spaces the limit is not necessarily round and even not necessarily Riemannian. We start from the convergence of the rescaled induced metric.

\begin{Theorem}
For any given $\kkk$ there is a smooth $\sss^a$-invariant function $f:\sss^m\rightarrow\mathbb R$ such that the rescaled induced metric $\tilde g_t=|\mm_t|^{-\frac{2}{m+a}}g_t$ converges, as $t$ goes to infinity, to the metric $e^{2f}\sigma_{\kkk}$.
\end{Theorem}
\proof
For any time $t$ let $\tilde{\rho}(t)$ be the radius of a geodesic sphere $\mathcal{B}_t$ such that $|\mm_t|=|\mathcal{B}_t|$. The mean curvature of $\mathcal{B}_t$ is $\tilde H=\hat H(\tilde{\rho})$, hence
$$
\frac{d\tilde{\rho}}{dt}=\frac{1}{\psi(\tilde H)},
$$
then $\tilde{\rho}=\frac{t}{\psi(m+a)}+o(1)$ as $t\rightarrow\infty$. Consider the function $\tilde f(x,t)=\rho(x,t)-\tilde{\rho}(t)$. We claim that $\tilde f$ converges to a smooth function $\tilde f_{\infty}$. In fact as a consequence of Lemma \ref{higher} we know that for any $N\in\mathbb N$ $\nabla^N\tilde f=\nabla^N\rho$ is uniformly bounded. Moreover by Proposition \ref{grad exp} we have
\begin{eqnarray*}
\left|\frac{\partial\tilde f}{\partial t}\right| & = & \left|\frac{v}{\psi(H)}-\frac{1}{\psi(\tilde H)}\right|\\
&	\leq & \frac{1}{\psi(H)}|v-1|+\frac{1}{\psi(H)\psi(\tilde H)}\left(\left|\psi(H)-\psi(m+a)\right|+\left|\psi(\tilde H)-\psi(m+a)\right|\right)\\
&\leq & ce^{-\frac{2}{\psi(m+a)}t}.
\end{eqnarray*}
When $\kkk=\rrr$ $e=\sigma_{\rrr}$ by Notation \ref{nnn}. When $\kkk\neq\rrr$ $e$ converges to $\sigma_{\kkk}$. Moreover by definition of $\tilde{\rho}$ we have that $|\mm_t|=\omega_m\sinh^m(\tilde{\rho})\cosh^a(\tilde{\rho})$ for some constant $\omega_m$. Therefore there is a constant $c>0$ such that when $t\rightarrow\infty$ we have
$$
\sinh(\rho)|\mm_t|^{-\frac{1}{m+a}}\approx ce^{\rho-\tilde{\rho}}=ce^{\tilde f}.
$$
It follows that there exists a positive constant $c$ such that
$$
\lim_{t\rightarrow\infty}\tilde g_t=ce^{2\tilde f_{\infty}}\sigma_{\kkk},
$$
from which it is easy to find the function $f$ that we are looking for.
\cvd

We want to describe the construction of examples such that the associated limit $\tilde g_{\infty}$ has not constant scalar curvature. The proof is an adaptation to a general $\psi$ of the techniques developed for the inverse mean curvature flow. As previous literature suggests, we need different tools for different values of $\kkk$. When $\kkk=\rrr$ we use the modified Hawking mass taken from \cite{HW}, but it is not useful in the other two ambient manifolds. In the other cases we use the Brown-York like masses introduced in \cite{Pi3,Pi4}.

\subsection{The case of the real hyperbolic space}
In this subsection we will focus only on the real hyperbolic space. Following Hung and Wang \cite{HW}, for any hypersurface $\mm$ in $\rrr\hhh^n$ we consider the modified Hawking mass
\begin{equation}\label{def HM}
Q(\mm):=|\mm|^{-1+\frac{4}{m}}\int_{\mm}\aao d\mu.
\end{equation}

\noindent With our notation Proposition 5 of \cite{HW} can be rewritten as follow.

\begin{Proposition}\label{HM lim}
\cite{HW} Let $f:\sss^{m}\rightarrow\mathbb R$ be a smooth function. For any $\tau$ big enough let $\widetilde{\mm}^{\tau}$ be the star-shaped hypersurface of $\mathbb{RH}^n$ defined by the radial function $\tilde{\rho}(z)=\tau+f(z)+o(1)$. Then
$$
\lim_{\tau\rightarrow\infty}Q(\widetilde{\mm}^{\tau})=\left(\int_{\sss^m}e^{mf}d\sigma\right)^{-1+\frac{4}{m}}\int_{\sss^m}e^{(m-2)f}\left|\accentset{\circ}{\nabla}^2e^{-f}\right|^2d\sigma.
$$
The limit of the rescaled induced metric is $e^{2f}\sigma$. It has constant scalar curvature if and only if $\lim_{\tau\rightarrow\infty}Q(\widetilde{\mm}^{\tau})=0$.
\end{Proposition}

We can compute the evolution of the modified Hawking mass under the flow \eqref{def flow}.

\begin{Lemma}\label{HM ev}
Let $\mm_t$ be a closed hypersurface of $\mathbb{RH}^n$ evolving according to \eqref{def flow}, then
\begin{eqnarray*}
|\mm_t|^{1-\frac{4}{m}}\frac{d Q(\mm_t)}{d t} & = & \left(\frac 4m-1\right)\left[\left|\mm_t\right|^{-1}\int_{\mm_t}\aao d\mu_t\int_{\mm_t}\left(\frac{H}{\psi}\right)d\mu_t-\int_{\mm_t}\aao\frac{H}{\psi}d\mu_t\right]\\
&&-\int_{\mm_t}\frac{2}{\psi}\accentset{\circ}{h}_i^j\accentset{\circ}{h}_j^k\accentset{\circ}{h}_k^i d\mu_t-2\int_{\mm_t}\frac{\psi'}{\psi^2}\nabla_iH\nabla^j\accentset{\circ}{h}_j^id\mu_t.
\end{eqnarray*}
\end{Lemma}
\proof
By Lemma \ref{eq evoluz} and the explicit expression of the curvature tensor of the ambient space we have 
\begin{eqnarray*}
\frac{\partial \accentset{\circ}{h}_i^j}{\partial t}&=&-\nabla_i\nabla^j\frac{1}{\psi}-\frac{1}{\psi}\left(\accentset{\circ}{h}_i^k\accentset{\circ}{h}_k^j+\frac{2H}{m}\accentset{\circ}{h}_i^j\right)+\left(\frac{1}{\psi}+\frac{\partial H}{\partial t}\right)\delta_i^j.
\end{eqnarray*}
Therefore
$$
\frac{\partial\aao}{\partial t}=2\nabla^j\left(\frac{\psi'\nabla_i H}{\psi^2}\right)\accentset{\circ}{h}_j^i-\frac{4H}{m\psi}\aao-\frac{2}{\psi}\accentset{\circ}{h}_i^j\accentset{\circ}{h}_j^k\accentset{\circ}{h}_k^i.
$$
The result follows easily considering the evolution of the volume form in Lemma \ref{eq evoluz}.
\cvd

The goal is to show that if $Q$ is decreasing, then it does so very slowly. Looking at the evolution of $Q$, we need to add an estimate to those of the previous Sections.

\begin{Lemma}\label{gtf}
If $\kkk=\rrr$, there is a positive constant $c$ such that
$$
\gtf\leq ce^{-\frac{6}{\pmm}t}.
$$
\end{Lemma}
\proof In this proof the $C_i$ will be positive constants. By Lemma \ref{eq evoluz}, Proposition \ref{grad exp} and Corollary \ref{conv 2ff} we can compute
\begin{eqnarray*}
\frac{\partial \aao}{\partial t}&=&\frac{\psi'}{\psi^2}\Delta\aao-2\frac{\psi'}{\psi^2}\gtf+2\frac{\psi''\psi-2(\psi')^2}{\psi^3}\nabla_iH\nabla_jH\tf_j^i\\
&&-2\left(\frac{1}{\psi}+\frac{H\psi'}{\psi^2}\right)\tf_i^j\tf_j^k\tf_k^i+2\left(\frac{\psi'}{\psi^2}\left(|A|^2+m-\frac 2m H\right)-\frac{2H}{m\psi}\right)\aao\\
&\leq&\frac{\psi'}{\psi^2}\Delta\aao-C_1\gtf+\left(C_2e^{-\frac{2}{\pmm}t}-\frac{4}{\pmm}\right)\aao+C_3e^{-\frac{6}{\pmm}t}.
\end{eqnarray*}
Moreover
\begin{eqnarray*}
\frac{\partial\gtf}{\partial t} & = & 2\nabla^s\tf_j^i\nabla_s\frac{\partial\tf_i^j}{\partial t}+\frac{\partial g^{rs}}{\partial t}\nabla_s\tf_i^j\nabla_r\tf_j^i\\
&=&2\nabla^s\tf_j^i\nabla_s\left[\frac{\psi'}{\psi^2}\Delta\tf_i^j+\frac{\psi''\psi-2(\psi')^2}{\psi^3}\nabla_iH\nabla^jH-\left(\frac{1}{\psi}+\frac{H\psi'}{\psi^2}\right)\tf_i^r\tf_r^j\right.\\
&&\phantom{2\nabla^s\tf_j^i\nabla_s}\left.+\left(\frac{\psi'}{\psi^2}\left(|A|^2+m-\frac{2H^2}{m\psi}\right)-\frac{2H}{m\psi}\right)\tf_i^j\right]-\frac{2}{\psi}h^{rs}\nabla_s\tf_i^j\nabla_r\tf_j^i\\
&\leq&\frac{\psi'}{\psi^2}\Delta\aao+\left(C_4e^{-\frac{2}{\pmm}t}-\frac{6}{\pmm}\right)\gtf+C_5\gtfq.
\end{eqnarray*}

Consider the auxiliary function $\beta=\log\gtf+K\aao$ for some positive constant $K$ to be determined later. From the two evolution equations above we get
$$
\frac{\partial\beta}{\partial t}\leq \frac{\psi'}{\psi^2}\Delta\beta+\frac{\psi'}{\psi^2}\left|\nabla\log\gtf\right|^2-\frac{6}{\pmm}+C_5\gtf-KC_1\gtf+C_6e^{-\frac{2}{\pmm}t}.
$$
By Corollary \ref{conv 2ff} there is a constant $C_7$ such that for any time $t$ we have $\aao\leq C_7e^{-\frac{4}{\pmm}t}$. Fix $t^*$ big enough such that
$$
16C_5C_7\frac{\psi'}{\psi^2}e^{-\frac{4}{\pmm}t^*}\leq C_1^2.
$$
Consider $\beta$ only for times $t\geq t^*$. In a point $(x_0,t_0)$ where $\beta$ attains its maximum we have
$$
\left|\nabla\log\gtf\right|^2=K^2\left|\nabla\aao\right|^2\leq 4K^2\aao\gtf\leq 4K^2C_7e^{-\frac{4}{\pmm}t^*}\gtf.
$$
By the choice of $t^*$ we can find a $K>0$ such that 
$$
4K^2C_7\frac{\psi'}{\psi^2}e^{-\frac{4}{\pmm}t^*}-KC_1+C_5\leq 0.
$$
By the maximum principle we have $\beta\leq-\frac{6}{\pmm}t+C_8$ and, by definition of $\beta$ the desired result follows.
\cvd

\begin{Proposition}\label{stima HM}
Let $\mm_t$ be the evolution of a star-shaped, mean convex hypersurface of $\mathbb{RH}^n$, then there is a positive constant $c$ such that
$$
\frac{d Q(\mm_t)}{dt}\geq -c e^{-\frac{2}{\psi(m)}t}.
$$
\end{Proposition}
\proof
By Corollary \ref{conv 2ff} we know that in the real hyperbolic space $\aao\leq ce^{-\frac{4}{\psi(m)}t}$. Moreover, by Corollary \ref{grad exp} we have
$$
\left||\mm|^{-1}\int\frac{H}{\psi}d\mu-\frac{H}{\psi}\right|\leq\left||\mm|^{-1}\int\frac{H}{\psi}d\mu-\frac{m}{\psi(m)}\right|+\left|\frac{m}{\psi(m)}-\frac{H}{\psi}\right|\leq ce^{-\frac{2}{\psi(m)}t}.
$$
Classical inequalities say that there is a positive constant $c$ such that
$$
|\nabla H|^2\leq m|\nabla A|^2\leq c\gtf.
$$
The first inequality is trivial, the second one follows by Lemma 2.2 of \cite{Hu}. Finally, we can use the Cauchy-Schwarz inequality and Lemma \ref{gtf} to estimate the last term in the evolution equation of Q:
$$
\left|\nabla_iH\nabla^j\tf_j^i\right|\leq |\nabla H|\cdot |\nabla\accentset{\circ}{A}|\leq c\gtf\leq ce^{-\frac{6}{\pmm}t}.
$$
\cvd

Now we have all the ingredients to repeat the construction of the not round examples described in \cite{HW}. We recall it briefly for completeness. Pick a constant $c_0>0$ big enough and a function $\bar f$ such that 
$$
\left(\int_{\sss^m}e^{m\bar f}d\sigma\right)^{-1+\frac{4}{m}}\int_{\sss^m}e^{(m-2)\bar f}\left|\accentset{\circ}{\nabla}^2e^{-\bar f}\right|^2d\sigma > 4c_0.
$$
Choose $\tau$ big enough such that the hypersurface $\widetilde{M}^{\tau}$ defined in Proposition \ref{HM lim} is mean convex and $Q(\widetilde{M}^{\tau})>2c_0$. Let $\mm^{\tau}_t $ be its evolution according to the flow \eqref{def flow}. We proved that this evolution is defined for any positive time and, as $t$ diverges, the rescaled induced metric converges to $e^{2f}\sigma$ for some function $f:\sss^m\rightarrow\mathbb R$. By Proposition \ref{stima HM} we have that $\lim_{t\rightarrow\infty}Q(\mm^{\tau}_t)>c_0>0$ and by Proposition \ref{HM lim} we can conclude that $e^{2f}\sigma$ is not round.

\subsection{The case of the complex and quaternionic hyperbolic space}
The richer geometry of $\ccc\hhh^n$ and $\hhh\hhh^n$ makes the research of a not round limit harder. In fact it is well known that in these spaces there are no totally umbilical hypersurfaces (see Theorem 5.1 of \cite{NR} for a proof). Therefore $\aao$ is always bounded away from zero and the modified Hawking mass is no more useful in this context. On the other hand the complexity of the problem can be reduced using the $\sss^a$-invariance. In fact, under this further hypothesis, Lemma \ref{Yamabe} suggests that the limits with constant scalar curvature should be very rare. In order to overcome all these difficulties in \cite{Pi3,Pi4} we defined the following weaker notion of mass.

\begin{equation}\label{def Q}
Q(\mm)=|\mm|^{-1+\frac{2}{m+a}}\int_{\mm}\left(H-\hat H\right)d\mu,
\end{equation}
where $\hat H$ is the function defined in \eqref{H hat}. Note that $Q$ in general does not have a sign, but it is bounded. Moreover it is important to keep in mind that, even if we are using an unique symbol, $Q$ depends on the choice of $\kkk$.

\begin{Proposition}\label{Q lim}
Let $f:\sss^{m}\rightarrow\mathbb R$ be a $\sss^a$-invariant smooth function. For any $\tau$ big enough let $\widetilde{\mm}^{\tau}$ be the star-shaped hypersurface of $\kh^n$ defined by the radial function $\tilde{\rho}(z)=\tau+f(z)+o(1)$. Then
$$
\lim_{\tau\rightarrow\infty}Q(\widetilde{\mm}^{\tau})=\left(\int_{\sss^m}e^{(m+a)f}d\sigma\right)^{-1+\frac{2}{m+a}}\int_{\sss^m}e^{(m+a)f}\left(e^{-f}\Delta e^{-f}-\frac{m+a}{2}|\nabla e^{-f}|^2\right)d\sigma.
$$
The limit of the rescaled induced metric on $\widetilde{\mm}^{\tau}$ is $e^{2f}\sigma_{\kkk}$. If $\lim_{\tau\rightarrow\infty}Q(\widetilde{\mm})\neq 0$ then $e^{2f}\sigma_{\kkk}$ does not have constant (Webster or qc) scalar curvature.
\end{Proposition}
We omit the proof of this result because it is a straightforward generalization of Proposition 9.1 of \cite{Pi4}: it does not depend on $\psi$ and only few minor changes are needed in order to include the quaternionic hyperbolic space. Now we need the evolution equation of $Q$.

\begin{Lemma}\label{Q ev}
Let $\mm_t$ be the evolution of a star-shaped, mean convex, $\sss^a$-invariant hypersurface of $\kh^n$, then
\begin{eqnarray*}
|\mm_t|^{1-\frac{2}{m+a}}\frac{d Q(\mm_t)}{d t} & = & \left(-1+\frac{2}{m+a}\right)|\mm_t|^{-1}\int_{\mm_t}\frac{H}{\psi}d\mu\int_{\mm_t}(H-\hat H)d\mu\\
&&+\int_{\mm_t}\left[\frac{H}{\psi}(H-\hat H)-\frac{1}{\psi}\left(|A|^2+\bar Ric(\nu,\nu)\right)\right]d\mu\\
&&+\int_{\mm_t}\frac{v}{\psi}\left(\frac{m}{\siq}-\frac{a}{\coq}\right)d\mu.
\end{eqnarray*}
\end{Lemma}
\proof
The result follows combining Lemma \ref{eq evoluz} with
$$
\frac{\partial\hat H}{\partial t} = \frac{v}{\psi}\left(\frac{a}{\coq}-\frac{m}{\siq}\right),
$$
and the fact that for any $t$
$$
\int_{\mm_t}\left(\frac{\psi'}{\psi^2}\Delta H+\frac{\psi''\psi-2(\psi')^2}{\psi^3}|\nabla H|^2\right)d\mu=-\int_{\mm_t}\Delta\left(\frac{1}{\psi}\right)d\mu=0.
$$
\cvd

With the help of the results of the previous Sections we can estimate the evolution equation of $Q$ proving that if $Q$ decays, then it does so very slowly.

\begin{Proposition}\label{stima Q}
Let $\mm_t$ be the evolution of a star-shaped, mean convex, $\sss^a$-invariant hypersurface of $\kh^n$, then there is a positive constant $c$ such that
$$
\frac{d Q(\mm_t)}{dt}\geq -c e^{-\frac{2}{\psi(m+a)}t}.
$$
\end{Proposition}
\proof
By the explicit expressions of the second fundamental form \eqref{2ff}, of the mean curvature \eqref{H} and by Lemma \ref{hessiani} we have
\begin{eqnarray*}
|A|^2+\bar Ric(\nu,\nu) & = & \frac{\hat{\varphi}_{ik}\hat{\varphi}_{js}\tilde{e}^{kj}\tilde{e}^{si}}{v^2\siq}+\frac{2\coh}{v\sih}\left(H-\frac{\hat H}{v}\right)\\
&&+m\left(\frac{\coq}{v^2\siq}-1\right)+a\left(\frac{\siq}{v^2\coq}+\frac{2}{v^2}-3\right)\\
& = & \frac{\varphi_{ik}\varphi_{js}\tilde{\sigma}^{kj}\tilde{\sigma}^{si}}{v^2\siq}+\frac{2\coh\hat H}{v^2(v+1)\sih}|\nabla\varphi|^2-\frac{m+a}{v^2}|\nabla\varphi|^2\\
&&+\frac{2\coh}{v\sih}\left(H-{\hat H}\right)+ \frac{m}{v^2\siq}-\frac{a}{v^2\coq}.
\end{eqnarray*}
By substituting this formula in the result of Lemma \ref{Q ev} we can rearrange the terms of the evolution equation of $Q$ as follows:
\begin{equation}\label{QQ}
\begin{array}{rcl}
|\mm_t|^{1-\frac{2}{m+a}}\frac{d Q(\mm_t)}{dt} & = &2\int\left(H-\hat H\right)\left(\frac{1}{(m+a)|\mm_t|}\int\frac{H}{\psi}d\mu-\frac{\coh}{v\psi\sih}\right)d\mu\\
&&+\int\left(H-\hat H\right)\left(\frac{H}{\psi}-\frac{1}{|\mm_t|}\int\frac{H}{\psi}d\mu\right)d\mu\\
&&+\int\frac{1}{\psi}\left(v-\frac{1}{v^2}\right)\left(\frac{m}{\siq}-\frac{a}{\coq}\right)d\mu\\
&&+\int\frac{|\nabla\varphi|^2}{v^2\psi}\left(m+a-\frac{2\coh\hat H}{(v+1)\sih}\right)d\mu-\int\frac{\varphi_{ik}\varphi_{js}\tilde{\sigma}^{ij}\tilde{\sigma}^{ks}}{v^2\psi\siq}d\mu.
\end{array}
\end{equation}
We claim that every term in the right hand side of \eqref{QQ} is smaller than $ce^{-\deo t}|\mm_t|$, for some constant $c>0$. In fact, by Proposition \ref{grad exp}, Lemma \ref{higher}, the facts that $\rho$ grows like $\frac{t}{\pma}$, and $\frac{1}{\psi}$ is bounded, the various terms of \eqref{QQ} can be estimate as follow:

$$
\begin{array}{rcl}
|H-\hat H|&\leq& |H-m-a|+|m+a-\hat H|\leq ce^{-\eo t};\\
&&\\
\left|\frac{1}{(m+a)|\mm_t|}\int\frac{H}{\psi}d\mu\right.-\left.\frac{\coh}{v\psi\sih}\right|&\leq&\frac{1}{(m+a)|\mm_t|}\int\left|\frac{H}{\psi}-\frac{m+a}{\psi(m+a)}\right|d\mu \\
&&+\frac{\coh}{v\sih}\left|\frac{1}{\psi}-\frac{1}{\psi(m+a)}\right|+\frac{1}{v\psi(m+a)}\left|\frac{\coh}{\sih}-1\right|\\
&\leq & ce^{-\eo t};\\
&&\\
\left|\frac{H}{\psi}-\frac{1}{|\mm_t|}\int\frac{H}{\psi}d\mu\right|&\leq &\left|\frac{H}{\psi}-\frac{m+a}{\pma}\right|+\frac{1}{|\mm_t|}\int\left|\frac{H}{\psi}-\frac{m+a}{\pma}\right|d\mu\\
&\leq& ce^{-\eo t};\\
&&\\
v-\frac{1}{v^2}&=&\frac{|\nabla\varphi|^2}{v^2(v+1)}\leq ce^{-\eo t};
\end{array}
$$

$$
\begin{array}{rcl}
\left|m+a-\frac{2\coh\hat H}{(v+1)\sih}\right|&\leq& \left|m+a-\hat H\right|+\hat H\left|1-\frac{2\coh}{(1+v)\sih}\right|\\
&\leq& ce^{-\eo t};\\
&&\\
\left|\varphi_{ik}\varphi_{js}\tilde{\sigma}^{ij}\tilde{\sigma}^{ks}\right|&\leq& ce^{-\eo t}.
\end{array}
$$
\cvd

As in Proposition 9.4 of \cite{Pi4} we can construct many $\sss^a$-invariant examples $\mm_0$ such that $\lim_{t\rightarrow\infty}Q(\mm_t)>0$, then by Proposition \ref{Q lim} the sub-Riemannian limit metric $\tilde g_{\infty}$ cannot have constant (Webster or qc) scalar curvature. The strategy is analogous to that of Hung and Wang \cite{HW} described in the previous subsection, therefore we omit it. It is sufficient to use the mass \eqref{def Q} instead of the modified Hawking mass \eqref{def HM}.

\bigskip

\noindent \textsc{Giuseppe Pipoli},\\ Department of Information Engineering, Computer Science and Mathematics, Universit\`a degli Studi dell'Aquila, via Vetoio 1, 67100 L'Aquila, Italy.\\ E-mail: giuseppe.pipoli@univaq.it

\end{document}